\documentclass{amsart}

\usepackage{fancyhdr}
\usepackage{lastpage}
\usepackage{stmaryrd,yhmath}

\newcommand{\bdis}{\begin{displaymath}}
\newcommand{\edis}{\end{displaymath}}
\newcommand{\be}{\begin{equation}}
\newcommand{\ee}{\end{equation}}
\newcommand{\mbb}{\mathbb}
\newcommand{\mcal}{\mathcal}

\DeclareMathOperator{\li}{li}

\theoremstyle{definition}

\newtheorem{cor}[]{Corollary}

\theoremstyle{remark}
\newtheorem{remark}[]{Remark}

\newtheorem*{mydef1}{{\bf Theorem}}

\numberwithin{equation}{section}

\begin{document}

\title{Riemann hypothesis and some new integrals connected with the integral negativity of the remainder in the formula for the prime-counting function $\pi(x)$}

\author{Jan Moser}

\address{Department of Mathematical Analysis and Numerical Mathematics, Comenius University, Mlynska Dolina M105, 842 48 Bratislava, SLOVAKIA}

\email{jan.mozer@fmph.uniba.sk}


\begin{abstract}
In this paper a new integral for the remainder of $\pi(x)$ is obtained. It is proved that there is an infinite set of the formulae containing miscellaneous parts of this
integral.
\end{abstract}

\maketitle

\section{The result}

Let us remind that

\bdis
\pi(x)=\int_0^x \frac{{\rm d}t}{\ln t}+P(x)=\li(x)+P(x)
\edis

where $\pi(x)$ is the prime-counting function. In 1914 Littlewood proved that the remainder $P(x),\ x\geq 2$ changes the sign infinitely many times
(see \cite{1}). However, on the Riemann hypothesis, the inequality

\be \label{1.1}
\int_2^Y P(x){\rm d}x<0 , \ Y\to\infty
\ee

holds true (see \cite{2}, p. 106), i.e. the remainder $P(x)$ is \emph{negative in mean}. In this direction the following theorem holds true.

\begin{mydef1}
Let $0<\Delta$ is a sufficiently small fixed number. Then, on Riemann hypothesis, there is a positive function $N(\delta),\ \delta\in (0,\Delta)$ such that
\be \label{1.2}
\int_2^{N(\delta)}P(x){\rm d}x=-\frac{a}{\delta}+\mcal{O}(1),\ \delta\in (0,\Delta),\ a=e^{9/2}
\ee
where $\mcal{O}(1)$ stands for a continuous and bounded function.
\end{mydef1}

\begin{remark}
The formula (\ref{1.2}) is the new asymptotic formula in the direction of (\ref{1.1}).
\end{remark}

\begin{remark}
Since (see (\ref{1.2})
\bdis
\liminf_{\delta\to 0^+} \int_2^{N(\delta)}P(x){\rm d}x=-\infty \ \Rightarrow \ \liminf_{\delta\to 0^+} N(\delta)=+\infty
\edis
then
\be \label{1.3}
\lim_{\delta\to 0^+}N(\delta)=\infty .
\ee
\end{remark}

\section{Some properties of asymptotic multiplicability for the integral (\ref{1.2})}

Let $0< \delta_k,\ k=1,\dots ,l$. Since (see (\ref{1.2}))
\bdis
\int_2^{N(\delta)}P(x){\rm d}x=-\frac{a}{\delta_k}\left\{ 1+o(1)\right\},\ \delta_k\to 0^+
\edis
then we have

\begin{cor}
If $\delta_k\to 0^+,\ k=1,\dots ,l$ then
\be \label{2.1}
\prod_{k=1}^l\int_2^{N(\delta_k)}P(x){\rm d}x\sim (-a)^{l-1}\int_2^{N\left(\prod_{k=1}^l\delta_k\right)}P(x){\rm d}x ,
\ee
i.e.
\be \label{2.2}
\prod_{k=1}^l\int_2^{N(\delta_k)}\frac{P(x)}{a}{\rm d}x\sim (-1)^{l-1}\int_2^{N\left(\prod_{k=1}^l\delta_k\right)}\frac{P(x)}{a}{\rm d}x ,
\ee
for arbitrary fixed $l\in\mbb{N}$.
\end{cor}

\begin{remark}
Let, for example,
\bdis
n=p_1^{\alpha_1}p_2^{\alpha_2}\dots p_l^{\alpha_l},\ p_k^{\alpha_k}\to\infty,\ k=1,\dots ,l
\edis
be the canonical product for $n$. Then (see (\ref{2.2}))
\be \label{2.3}
\prod_{k=1}^l\int_2^{N(p_k^{-\alpha_k})}\frac{P(x)}{a}{\rm d}x\sim (-1)^{l-1}\int_2^{N\left(\prod_{k=1}^l p_k^{-\alpha_k}\right)}\frac{P(x)}{a}{\rm d}x .
\ee
\end{remark}

\section{The division in the set of some integrals}

The simple algebraic formula

\bdis
\left(\frac{1}{\delta_2}\right)^n-\left(\frac{1}{\delta_1}\right)^n=\left(\frac{1}{\delta_2}-\frac{1}{\delta_2}\right)
\sum_{k=0}^{n-1}\left(\frac{1}{\delta_1}\right)^k\left(\frac{1}{\delta_2}\right)^{n-k-1},\ \delta_2<\delta_1 ,
\edis
implies due to (\ref{1.2}) the following dual formula.

\begin{cor}
\be\label{3.1}
\frac{\int_{N(\delta_1^n)}^{N(\delta_2^n)}\frac{P(x)}{a}{\rm d}x}{\int_{N(\delta_1)}^{N(\delta_2)}\frac{P(x)}{a}{\rm d}x}\sim
-\sum_{k=0}^{n-1}\int_2^{N(\delta_1^k\delta_2^{n-k-1})}\frac{P(x)}{a}{\rm d}x ,
\ee
for arbitrary fixed $n\in\mbb{N},\ \delta_1,\delta_2\to 0^+$ .
\end{cor}

\begin{remark}
The asymptotic division
\bdis
\frac{\int_{N(\delta_1^n)}^{N(\delta_2^n)}\frac{P(x)}{a}{\rm d}x}{\int_{N(\delta_1)}^{N(\delta_2)}\frac{P(x)}{a}{\rm d}x}
\edis
is defined by the formula (\ref{3.1}).
\end{remark}

\section{The division of some integrals on asymptotically equal parts}

First of all we have (see (\ref{1.2}))
\be \label{4.1}
\int_2^{N\left(\frac{a\delta}{n}\right)}P(x){\rm d}x\sim -\frac{n}{\delta},\ n\in\mbb{N},\ \delta\in (0,\Delta) .
\ee
Then from (\ref{4.1}) the following corollary follows.

\begin{cor}
\be \label{4.2}
\int_{N\left(\frac{a\delta}{k}\right)}^{N\left(\frac{a\delta}{k+1}\right)}P(x){\rm d}x\sim
\int_{N\left(\frac{a\delta}{l}\right)}^{N\left(\frac{a\delta}{l+1}\right)}P(x){\rm d}x ,
\ee
for every $k,l\in\mbb{N},\ \delta\to 0^+$.
\end{cor}

\begin{remark}
Thus the sequence
\bdis
\left\{N\left(\frac{a\delta}{n}\right)\right\}_{n=1}^\infty
\edis
subdivides each of the integrals
\bdis
\int_{N\left(\frac{a\delta}{n}\right)}^{N\left(\frac{a\delta}{n+m}\right)}P(x){\rm d}x,\ \forall \ n,m\in\mbb{N}
\edis
on the asymptotically equal parts (see (\ref{4.2})) .
\end{remark}

\begin{remark}
For every sequence of the type
\bdis
\delta_k(\alpha_0,r)=\frac{a\delta}{\alpha_0+kr},\ k=0,1,2,\dots ,\ \alpha_0,r>0
\edis
we obtain the similar result.
\end{remark}

\begin{remark}
Thus we see that there is an infinite set of the formulae containing miscellaneous parts of the integral
\bdis
\int_2^{N(\delta)}P(x){\rm d}x
\edis
(see for example (\ref{1.3}), (\ref{2.1})-(\ref{2.3}), (\ref{4.2})).
\end{remark}

\section{Proof of the Theorem}

Let us remind the formula
\be \label{5.1}
\int_2^\infty\frac{\ln(e^{-2}x)}{x^{3/2+\delta}}P(x){\rm d}x=-\frac 1\delta+\mcal{O}(1),\ \delta\in (0,\Delta)
\ee
where $\mcal{O}(1),\ \delta\in (0,\Delta)$ is a continuous and bounded function (see \cite{4}). Using the von Koch estimate
$P(x)=\mcal{O}(\sqrt{x}\ln x)$, which follows from the Riemann hypothesis, we obtain (comp. \cite{4}, part 3)
\bdis
\left|\int_Y^\infty\frac{\ln(e^{-2}x)}{x^{3/2+\delta}}P(x){\rm d}x\right|< \frac A\delta Y^{-\delta/2}\leq 1
\edis
for
\bdis
Y\geq T(\delta)=\left(\frac A\delta\right)^{2/\delta},\ \delta\in (0,\Delta) .
\edis
Then we have (see (\ref{5.1}))
\be \label{5.2}
\int_{e^3}^Y\frac{\ln(e^{-2}x)}{x^{3/2+\delta}}P(x){\rm d}x=-\frac 1\delta+\mcal{O}(1),\ \delta\in (0,\Delta) ,
\ee
where $\int_2^{e^3} \ =\mcal{O}(1)$. Let us remind further the Bonet's form of the second mean-value theorem
\bdis
\int_a^b f(x)g(x){\rm d}x=f(a+0)\int_a^\xi g(x){\rm d}x ,\ \xi \in (a,b)
\edis
where $f(x),g(x)$ are integrable function on $[a,b]$, and $f(x)$ is a non-negative and non-increasing function (in our case $g(x)=P(x)$). We have
as a consequence the following formula (comp. \cite{3})
\be \label{5.3}
\int_{e^3}^Y\frac{\ln(e^{-2}x)}{x^{3/2+\delta}}P(x){\rm d}x=\frac{1}{e^{9/2+3\delta}}\int_{e^3}^{M(Y,\delta)} P(x){\rm d}x
\ee
where $M\in (e^3,Y)$. Let $Y=T(\delta)$, i.e.
\be \label{5.4}
M(Y,\delta)=M[T(\delta),\delta]=N(\delta) .
\ee
Then from (\ref{5.2}) by (\ref{5.3}) and (\ref{5.4}) the formula
\be \label{5.5}
\int_{e^3}^{N(\delta)} P(x){\rm d}x=-\frac{e^{9/2+3\delta}}{\delta}+\mcal{O}(1),\ \delta\in (0,\Delta)
\ee
follows. Since
\bdis
e^{3\delta}=1+\mcal{O}(\delta),\ \int_2^{e^3}P(x){\rm d}x=\mcal{O}(1),
\edis
we obtain our formula (\ref{1.2}) from (\ref{5.5}).

\section{Concluding remarks: on the Littlewood's point and on the point of simple discontinuity of the function $N(\delta)$}

Since the right-hand side of the formula (\ref{5.5}) contains the continuous function ($Y=T(\delta)$) then
\bdis
\int_{T(\delta)}^\infty \frac{\ln (e^{-2}x)}{x^{3/2+\delta}}P(x){\rm d}x
\edis
is also the continuous function. Let us remind that the continuity of the two functions $F(u)$ and $F(f(x))$ does not imply the
continuity of the function $f(x)$, (comp. $\sin\{\pi[x]\},\ x\in\mbb{R}$). \\

If $\delta_0\in (0,\Delta)$ is the point of the simple discontinuity of $N(\delta)$ (another kind of discontinuity of $N(\delta)$ is excluded), i.e. if
$N(\delta_0+0)\not= N(\delta_0-0)$ then we have (see (\ref{1.2}))
\bdis
\int_{N(\delta_0-0)}^{N(\delta_0+0)} P(x){\rm d}x=0 ,
\edis
and from this
\be \label{6.1}
\mu[P(x)]\left\{N(\delta_0+0)-N(\delta_0-0)\right\}=0 \ \Rightarrow \ \mu[P(x)]=0
\ee
follows, where $\mu$ is the mean-value of $P(x)$ relative to the segment generated by the points $N(\delta_0-0),N(\delta_0+0)$. \\

Let us remind that the remainder $P(x)$ changes its sign at the Littlewood' sequence points $\{ L_k\}_{k=1}^\infty$.

\begin{remark}
If $\delta_0\in (0,\Delta)$ is the point of the simple discontinuity of $N(\delta)$ (if any), then from (\ref{6.1}) we have: there is $L_k$ for which
\bdis
L_k\in (N(\delta_0-0),N(\delta_0+0)),\ [N(\delta_0-0),N(\delta_0+0)]
\edis
holds true.
\end{remark}

\end{document}